\documentclass[11pt]{article}
\usepackage{amsthm,amstext}
\usepackage{amsmath}
\usepackage{graphicx}
\usepackage{amsfonts}
\usepackage{amssymb}
\usepackage[colorlinks,linkcolor=blue,citecolor=blue]{hyperref}
\usepackage{tikz}
\usepackage{enumerate}
\usepackage{mathrsfs}
\usepackage{rotate}
\usepackage{fancyhdr}

\theoremstyle{plain}
\newtheorem{theorem}{Theorem}[section]
\newtheorem{lemma}[theorem]{Lemma}
\newtheorem{proposition}[theorem]{Proposition}

\newtheorem{corollary}[theorem]{Corollary}

\theoremstyle{definition}

\newcommand{\identity}{\mathrm{id}}

\newcommand{\rank}{\mathrm{rank~}}
\newcommand{\im}{\mathrm{im~}}

\newcommand{\collapse}{\mathrm{c}}
\newcommand{\nb}{\mathrm{nb}}
\newcommand{\Inj}{\mathrm{Inj}}
\newcommand{\Sing}{\mathrm{Sing}}

\date{}

\title{\bf Generating Sets of an Infinite Semigroup of Transformations Preserving a Zig-zag Order }\vspace{.25 in}
\begin{document}
\maketitle

\begin{center}
Laddawan Lohapan\\
Department of Mathematics, Faculty of Science\\
Khon Kaen University, Khon Kaen 40002, Thailand\\
email: lohapan$\_$l@kkumail.com
\end{center}
\begin{center}
J\"org Koppitz\\
Institute of Mathematics\\
Bulgarian Academy of Sciences\\
Acad. G. Bonchev St, Bl. 8, Sofia, 1113, Bulgaria\\
email: koppitz@math.bas.bg
\end{center}
\begin{center}
Somnuek Worawiset\\
Department of Mathematics, Faculty of Science\\
Khon Kaen University, Khon Kaen 40002, Thailand\\
email: wsomnu@kku.ac.th 
\end{center}

\begin{abstract}
A zig-zag order is like a directed path, only with alternating directions. A generating set of minimal size for the semigroup of all full transformations on a finite set preserving the zig-zag order was determined by Fenandes et al. in 2019. This paper deals with generating sets of the semigroup $F_{\mathbb{N}}$ of all full transformations on the set of all natural numbers preserving the zig-zag order. We prove that $F_{\mathbb{N}}$ has no minimal generating sets and present two particular infinite decreasing chains of generating sets of $F_{\mathbb{N}}.$
\end{abstract}
$Key~words~and~phrases:$ fence, zig-zag order, order-preserving, generating set, transformation
\section{Introduction}
This paper deals with generating sets of transformation semigroups. A full transformation on a set $X$ is a self-mapping on $X.$ The set of all full transformations on $X$ forms a semigroup $T_{X}$ under the usual composition of mappings. If $X$ is the $n$-element set $\{1,2,\ldots,n\},$ then we write $T_{n}$ rather than $T_{X}.$ In particular, $T_{n}$ is a finite semigroup of full transformations, which is the disjoint union of the symmetric group and the singular part $\Sing_{n}.$ In fact, $\Sing_{n}$ is an ideal of $T_{n}$ consisting of all full transformations with $\rank<n.$ The semigroup $\Sing_{n}$ is generated by the idempotents of $\rank n-1$ \cite{A}. Ayik et al. found a necessary and sufficient condition for any set of full transformations with $\rank n-1$ to be a generating set of $\Sing_{n}$ \cite{B}. The generating sets of the ideals  $K(n,r),r\in\{1,2,\ldots,n-1\},$ of $\Sing_{n}$ were determined by Ayik and Bugay \cite{C}. 

The set $O_{n}$ of all order-preserving full transformations on $\{1,2,\ldots,n\}$ with respect to the usual linear order on the natural numbers forms a semigroup, which is the disjoint union of the identity mapping on $\{1,2,\ldots,n\}$ and the singular part. The minimal size of a generating set of $O_{n}$ (i.e., the $\rank$ of $O_{n}$) is $n$ while the singular part is generated by its idempotents of $\rank n-1$ \cite{D}. A necessary and sufficient condition for any set of full transformations in the ideal $O(n,r), r\in\{1,2,\ldots,n-1\},$ to be a generating set of $O(n,r)$ was provided by Ayik and Bugay \cite{E}.
 
Generating sets for other (finite) semigroups of full transformations have been determined by several authors. Among these semigroups is the semigroup $F_{n}$ of all full transformations on $\{1,2,\ldots,n\}$ preserving the zig-zag order. Recall that the zig-zag order is a partial order, which is like a path, only with alternating directions. Full transformations on $\{1,2,\ldots,n\}$ preserving the zig-zag order were first studied by Currie and Visentine \cite{3} and Rutkowski \cite{16} in 1991 and 1992, respectively. In both papers, the authors calculated the cardinality of $F_{n},$ depending on the parity of $n.$ In \cite{F}, Fernandes, Koppitz, and Musunthia determined a generating set of $F_{n}$ of minimal size and gave a formula to calculate the rank of $F_{n}.$ Algebraic properties of $F_{n}$ were investigated by several authors in the last decade (e.g. \cite{I,H,G}).

Recall that uncountable semigroups have only uncountable generating sets. In order to make the situation more comfortable, Ru\v{s}kuc introduced the concept of a relative generating set (i.e. a relative rank) \cite{Firstgen}. For example, in \cite{J,K}, the authors considered the uncountable semigroup $T_{\mathbb{N}}$ and the semigroup $O_{\mathbb{N}}$ of all order-preserving full transformations on the set $\mathbb{N}$ of all natural numbers with respect to the usual linear order on $\mathbb{N}.$ One needs only one $\alpha\in T_{\mathbb{N}}\setminus O_{\mathbb{N}}$ such that $O_{\mathbb{N}}\cup \{\alpha\}$ generates $T_{\mathbb{N}},$ i.e. the relative rank of $T_{\mathbb{N}}$ modulo $O_{\mathbb{N}}$ is one, where $\{\alpha\}$ is said to be a relative generating set of $T_{\mathbb{N}}$ modulo $O_{\mathbb{N}}.$ On the other hand, in \cite{K}, Higgins, Mitchell, and Ru\v{s}kuc considered the set $C$ of all contractions on $\mathbb{N}$ and obtained that the relative rank of $T_{\mathbb{N}}$ modulo $C$ is uncountable. Also in \cite{K}, the authors pointed out that  the relative rank of $T_{\mathbb{N}}$ modulo a so-called dominated set is uncountable. 

In the present paper, we consider an extension of the zig-zag order on $\{1,2,\ldots,n\}$ to the set of all natural numbers $\mathbb{N}$. Let
\[
\begin{array}{ccc}
n\prec n+1 & \text{if } n\text{ is odd;} \\ 
n+1\prec n & \text{otherwise.}%
\end{array}%
\]%
The binary relation $\prec$ together with the diagonal on $\mathbb{N}$ is a partial order on $\mathbb{N},$ in fact, $\preceq$ is called the zig-zag order on $\mathbb{N}.$ Any element in the partially ordered set $(\mathbb{N},\preceq),$ which is called a fence, is either minimal or maximal. The set $F_{\mathbb{N}}$ of all full transformations on $\mathbb{N}$ preserving the zig-zag order forms a subsemigroup of $T_{\mathbb{N}}.$ Corollary 2.2. in \cite{K} and the fact that $F_{\mathbb{N}}$ is dominated imply that the relative rank of $T_{\mathbb{N}}$ modulo $F_{\mathbb{N}}$ is uncountable infinite. In fact, the study of the semigroup $F_{\mathbb{N}}$ extends the study of $F_{n}$ on another level (we have now an uncountable semigroup of full transformations). Furthermore, congruences on $F_{\mathbb{N}}$ were already determined in \cite{Lohapan}. Hence, a more detailed study of the semigroup $F_{\mathbb{N}}$ seems reasonably enough. An investigation of generating sets of $F_{n}$ will be provided in this paper.

Besides the zig-zag order $\preceq$ on $\mathbb{N},$ we also deal with the usual liner order $\leq$ on $\mathbb{N}.$ Excluding any confusion, we introduce the following agreements. Let $A$ be a non-empty subset of $\mathbb{N}.$ We use $\min (A)$ and $\max(A)$ for the smallest and the greatest element (if exists), respectively, in $A$ with respect to $\leq.$ Moreover, $A$ is said to be convex if $A$ is an interval with respect to $\leq.$ Note that the image of $\alpha$ (in symbols: $\im\alpha$) is a convex set. For $B\subseteq \mathbb{N},$ we write $A<B$ if $a<b$ for all $a\in A$ and all $b\in B.$

In the next section, we show that any transformation in $F_{\mathbb{N}}$ can be expressed as the product of one element from each of the sets $$\Theta:=\{\alpha\in F_{\mathbb{N}}:a\alpha^{-1}~\text{is a convex }\text{set for all}~a\in\im\alpha\}~\text{and}$$ $$\Lambda_{n}:=\{\alpha\in F_{\mathbb{N}}:\left\vert\nb(\alpha)\right\vert=0,~\collapse(\alpha)>0,\\1\alpha\geq n,~\text{and}~\left\vert \{1,2,\ldots,n\}\alpha\right\vert=n\}$$
for any $n\in \mathbb{N},$ where  
$$\nb(\alpha):=\{a\in\mathbb{N}:a\alpha=(a+1)\alpha\}~\text{and}$$ $$c(\alpha):=\left\vert\bigcup \left\{ a\alpha ^{-1}:a\in \im\alpha \text{ and }\left\vert a\alpha ^{-1}\right\vert \geq 2\right\}\right\vert.$$ 
Obviously, $c(\alpha)\leq c(\alpha\beta)$ for all $\alpha,\beta\in F_{\mathbb{N}}$ and $\collapse(\alpha)=0$ if and only if $\alpha$ is injective. It is worth mentioning that $F_{\mathbb{N}}$ has no minimal generating sets. The main paper's purpose is to give two particular infinite decreasing chains of generating sets of $F_{\mathbb{N}},$ which will be provided in Section 3.

Let $\alpha \in F_{\mathbb{N}}.$ The rank of $\alpha,$ (in symbols: $\rank\alpha)$ is the size of the image of $\alpha.$ Then $\rank\alpha$ can be finite (in symbols: $\rank\alpha <\aleph _{0}$) or countable
infinite (in symbols: $\rank\alpha =\aleph _{0}$). The set of all transformations in $F_{\mathbb{N}}$ with countable infinite rank will be
denoted by $F_{\mathbb{N}}^{\inf}$. For $n\in \mathbb{N}$, let $\Theta_{n}=\Theta\cap\Omega_{n,}$ where $$\Omega_{n}:=\{\alpha\in F_{\mathbb{N}}:1\alpha\geq n~\text{and}~\left\vert \{1,2,\ldots,n\}\alpha\right\vert=n\}.$$
Then we obtain that $\Lambda _{n}=\Lambda \cap \Omega _{n},$ where $\Lambda:=\{\alpha\in F_{\mathbb{N}}:\left\vert\nb(\alpha)\right\vert=0~\text{and}~\collapse(\alpha)>0\}.$ 
Just for convenience, for $\alpha\in F_{\mathbb{N}},$ we define the following sets, which will be used
subsequently:\\
$~$\\
$M_{\alpha}^{n}:=\{X\subseteq \mathbb{N}: \left\vert X\right\vert=n~\text{and}~X~\text{is a maximal convex set with respect to} \left\vert X\alpha\right\vert=1\};\\
M_{\alpha}:=\bigcup_{n\in\mathbb{N}}M_{\alpha}^{n};\\
M_{\alpha}^{*}:=M_{\alpha}\setminus M_{\alpha}^{1};\\
MS_{\alpha}^{n}:=\{X\subseteq \bigcup M_{\alpha}^{1}: X~\text{is a maximal convex set and} \left\vert X\right\vert=n\};\\
MS_{\alpha}:=\bigcup_{n\in\mathbb{N}}MS_{\alpha}^{n}.$\\

\noindent More in detail, a convex set $X\subseteq\mathbb{N}$ belongs to $M_{\alpha}^{n}$ if and only if $\left\vert X\right\vert=n,\left\vert X\alpha\right\vert=1,$ and $\left\vert Y\alpha\right\vert>1$ for any convex set $Y\subseteq\mathbb{N}$ with $X\subsetneq Y.$ Moreover, a convex set $X\subseteq\bigcup M_{\alpha}^{1}$ belongs to $SM_{\alpha}^{n}$ if and only if $\left\vert X\right\vert=n$ and $Y\nsubseteq\bigcup M_{\alpha}^{1}$ for any convex set $Y\subseteq\mathbb{N}$ with $X\subsetneq Y.$ For any $\beta\in F_{\mathbb{N}}$, it is clear that  $M_{\alpha}=M_{\beta}$ if and only if $M_{\alpha}^{*}=M_{\beta}^{*}.$\\
\noindent Further, let $C_{m}:=\{X: X\subseteq\{m,m+1,\ldots\}\}$ for all $m\in\mathbb{N}.$   
\section{On Minimal Generating Sets of $F_{\mathbb{N}}$}
First, we describe any transformation $\alpha$ in $F_{\mathbb{N}},$ that is, $\alpha$ preserves the partial order $\preceq$ on $\mathbb{N}.$ If $x,y\in\mathbb{N}$ with $x\prec y,$ then $x$ is odd and $y$ is even. Moreover, $x$ is the successor of $y$ or conversely $y$ is the successor of $x,$ which implies $\left\vert x-y\right\vert=1.$ When we apply $\alpha$ to both $x$ and $y,$ their images are related with respect to $\preceq,$ that is, $\left\vert x\alpha-y\alpha\right\vert\leq 1.$ This fact will be used subsequently without mentioning. Now, we characterize the elements of $F_{\mathbb{N}}$ by two properties, which are easy to verify. 
\begin{proposition}\label{FirstProposition}
Let $\alpha\in T_{\mathbb{N}}.$ Then $\alpha\in F_{\mathbb{N}}$ if and only if
\begin{itemize}
\item[$(i)$] $\left\vert x\alpha-(x+1)\alpha \right\vert\leq 1$ for all $x\in\mathbb{N};$
\item[$(ii)$] $x$ and $x\alpha$ have the same parity or $(x-1)\alpha=x\alpha=(x+1)\alpha$ for all $x\in\mathbb{N}\setminus\{1\}.$
\end{itemize}
\begin{proof}
Suppose $\alpha\in F_{\mathbb{N}}.$\\
(i) Let $x\in\mathbb{N}.$ Then $x\prec x+1$ or $x+1\prec x.$ Since $\alpha\in F_{\mathbb{N}},$ we obtain $x\alpha\preceq (x+1)\alpha$ and $(x+1)\alpha\preceq x\alpha,$ respectively. Then $\left\vert x\alpha-(x+1)\alpha\right\vert\leq 1.$\\ (ii) Suppose that there exists $x\in\mathbb{N}\setminus\{1\}$ such that $x$ and $x\alpha$ have different parities. Without loss of generality, suppose that $x$ is odd and $x\alpha$ is even. Assume $(x-1)\alpha\not=x\alpha.$ Then (i) implies $(x-1)\alpha\in\{x\alpha-1,x\alpha+1\}.$ It follows that $(x-1)\alpha$ is odd. This shows that $x\prec x-1$ but $(x-1)\alpha\prec x\alpha,$ that is, $\alpha\not\in F_{\mathbb{N}},$ a contradiction. Hence, $(x-1)\alpha=x\alpha.$ Similarly, we can show that $(x+1)\alpha=x\alpha.$\\

Conversely, suppose that (i) and (ii) hold. Let $x,y\in\mathbb{N}$ be such that $x\prec y.$ Then $x$ is odd and $y$ is even with $x\in\{y-1,y+1\}.$ By (i), we obtain $\left\vert x\alpha-y\alpha\right\vert\leq 1.$ It is enough to consider the case $\left\vert x\alpha-y\alpha\right\vert= 1.$ Since $x\in\{y-1,y+1\}$ and $\left\vert x\alpha-y\alpha\right\vert=1,$ we obtain that $y$ and $y\alpha$ are even by (ii) and so $x\alpha\prec y\alpha.$ Altogether, we conclude $x\alpha\preceq y\alpha.$ Therefore, $\alpha\in F_{\mathbb{N}}.$
\end{proof}
\end{proposition}
An immediate consequence of Proposition \ref{FirstProposition} is that $\left\vert A\right\vert$ is odd for all $A\in M_{\alpha}^{*}$ with $1\not\in A.$ In the following, we will use this fact as well as Proposition \ref{FirstProposition} without further mentioning. Any element in $F_{\mathbb{N}}$ can be described as the product of one element from each of the sets $\Theta$ and $\Lambda_{n}$ for any $n\in\mathbb{N}.$
\begin{proposition}\label{Proposition4}
$F_{\mathbb{N}}=\Theta\Lambda_{n}=\{\gamma_{1}\gamma_{2}:\gamma_{1}\in\Theta,\gamma_{2}\in\Lambda_{n}\}$ for all $n\in\mathbb{N}.$
\begin{proof}
Let $n\in\mathbb{N}$ and $\alpha\in F_{\mathbb{N}}.$ Then we consider the following two cases.\\
\textbf{Case 1:} $\left\vert M_{\alpha}\right\vert=\aleph_{0}.$ Suppose $M_{\alpha}=\{A_{i}:i\in\mathbb{N}\}$ with $A_{i}<A_{i+1}$ for all $i\in\mathbb{N}.$ Then $\left\vert A_{i}\right\vert<\aleph_{0}$ for all $i\in\mathbb{N}.$ For all $i\in\mathbb{N},$ let $m_{i}=\max(A_{i}).$ This means $A_{i}\alpha=\{m_{i}\alpha\}$ for all $i\in\mathbb{N}.$ Obviously, $\alpha\in F_{\mathbb{N}}$ and $\left\vert A_{i}\alpha\right\vert=1$ for all $i\in\mathbb{N}$ imply that for all $i\in\mathbb{N},$
\begin{equation}\label{equation1}
m_{i}~\text{and}~m_{i}\alpha~\text{have the same parity and}~\left\vert m_{i}\alpha-m_{i+1}\alpha\right\vert=1.
\end{equation}  
Let $k\in\mathbb{N}\setminus\{1,2,\ldots,n\}$ be such that $k$ and $m_{1}\alpha$ have the same parity. We define $\gamma_{1}:\mathbb{N}\to\mathbb{N}$ by
$$x\gamma_{1}:=k+i-1~\text{for all}~x\in A_{i}, i\in\mathbb{N}.$$
The transformation $\gamma_{1}$ is well defined since $\bigcup_{i\in\mathbb{N}}A_{i}=\mathbb{N}.$ Moreover, $A_{i}\gamma_{1}=\{k+i-1\}$ for all $i\in\mathbb{N}$ and thus, $M_{\gamma_{1}}=M_{\alpha}.$ It is clear that $\left\vert x\gamma_{1}-(x+1)\gamma_{1}\right\vert\leq 1$ for all $x\in\mathbb{N}.$ Since $k$ and $m_{1}\alpha$ have the same parity and $M_{\gamma_{1}}=M_{\alpha},$ we obtain that $x$ and $x\gamma_{1}$ have the same parity or $(x-1)\gamma_{1}=x\gamma_{1}=(x+1)\gamma_{1}$ for all $x\in\mathbb{N}\setminus\{1\}.$ Since $y\gamma_{1}^{-1}$ is a convex set for all $y\in\im\gamma_{1},$ we obtain $\gamma_{1}\in\Theta.$ Further, we define $\gamma_{2}:\mathbb{N}\to\mathbb{N}$ by
$$x\gamma_{2}:=\begin{cases}
m_{1}\alpha+k-x&\text{if}~x\in\{1,2,\ldots,k-1\};\\
m_{x-k+1}\alpha~~&\text{if}~x\in\{k,k+1,\ldots\}.
\end{cases}$$
By (\ref{equation1}) and the fact that $k$ and $m_{1}\alpha$ have the same parity, we can conclude that (i) and (ii) in Proposition \ref{FirstProposition} are satisfied for $\gamma_{2},$ that is, $\gamma_{2}\in F_{\mathbb{N}}.$ If $\rank\alpha=\aleph_{0},$ then there exists $y\in\{m_{2}\alpha,m_{3}\alpha,\ldots\}$ with $y=m_{1}\alpha+1,$ that is, $\gamma_{2}$ is not injective. If $\rank\alpha<\aleph_{0},$ then it is clear that $\gamma_{2}$ is not injective. Moreover, we have $\left\vert\nb(\gamma_{2})\right\vert=0, \left\vert\{1,2,\ldots,n\}\gamma_{2}\right\vert=n,$ and $1\gamma_{2}=m_{1}\alpha+k-1\geq k> n.$ Thus, $\gamma_{2}\in\Lambda_{n}.$ By straightforward calculations, we obtain $A_{i}\gamma_{1}\gamma_{2}=\{m_{i}\alpha\}$ for all $i\in\mathbb{N}.$ This shows $\gamma_{1}\gamma_{2}=\alpha.$\\
$~$\\
\textbf{Case 2:} $\left\vert M_{\alpha}\right\vert<\aleph_{0}.$ Suppose $M_{\alpha}=\{A_{i}:1\leq i\leq l\}$ for some $l\in\mathbb{N}$ with $A_{i}<A_{j}$ for all $1\leq i<j\leq l.$ Then $\left\vert A_{i}\right\vert<\aleph_{0}$ for all $i\in\mathbb{N}\setminus\{l,l+1,\ldots\}$ and $\left\vert A_{l}\right\vert=\aleph_{0}.$ Let $m_{i}=\max(A_{i})~\text{for all }i\in\mathbb{N}\setminus\{l,l+1,\ldots\}~\text{and}~m_{l}=\min(A_{l}).$
Then $A_{i}\alpha=\{m_{i}\alpha\}$ for all $i\in\{1,2,\ldots,l\}.$ Since $\alpha\in F_{\mathbb{N}}$ and $\left\vert A_{i}\alpha\right\vert=1$ for all $1\leq i\leq l,$ the following properties hold:
\begin{itemize}
\item[(a1)] $\left\vert m_{i}\alpha-m_{i+1}\alpha\right\vert=1$ for all $i\in\mathbb{N}\setminus\{l,l+1,\ldots\};$
\item[(a2)] $m_{i}$ and $m_{i}\alpha$ have the same parity for all $1\leq i\leq l,$ whenever $l>1.$
\end{itemize}
Let $k\in\mathbb{N}\setminus\{1,2,\ldots,n\}$ be such that $k$ and $m_{1}\alpha$ have the same parity. Then we define $\gamma_{1}:\mathbb{N}\to\mathbb{N}$ by $$x\gamma_{1}:=k+i-1~\text{for all}~x\in A_{i},1\leq i\leq l.$$
The transformation $\gamma_{1}$ is well defined since $\bigcup_{i\in \mathbb{N}}A_{i}=\mathbb{N}.$ Moreover, $A_{i}\gamma_{1}=\{k+i-1\}$ for all $1\leq i\leq l.$ Using the same arguments as in Case 1, we get $\gamma_{1}\in F_{\mathbb{N}}.$ Since $y\gamma_{1}^{-1}$ is a convex set for all $y\in\im\gamma_{1},$ we have $\gamma_{1}\in\Theta.$ 
Further, let $\gamma_{2}:\mathbb{N}\to\mathbb{N}$ by
$$x\gamma_{2}:=\begin{cases}
m_{1}\alpha+k-x&\text{if }x\in\{1,2,\ldots,k-1\};\\
m_{x-k+1}\alpha~~&\text{if }x\in\{k,k+1,\ldots,k+l-1\};\\
m_{l}\alpha+x-k-l+1&\text{if }x\in\{k+l,k+l+1,\ldots\}.
\end{cases}$$
By (a1), we have $\left\vert x\gamma_{2}-(x+1)\gamma_{2}\right\vert\leq 1$ for all $x\in\mathbb{N}.$ Moreover, $x$ and $x\gamma_{2}$ have the same parity for all $x\in\mathbb{N}$ by (a2) and the property of $k.$ Hence, $\gamma_{2}\in F_{\mathbb{N}}.$ Since $\im\gamma_{2}=\{m_{1}\alpha,\ldots,m_{l}\alpha,m_{l}\alpha+1,m_{l}\alpha+2,\ldots\}$ is a convex set, $\rank\gamma_{2}=\aleph_{0},$ and $k\gamma_{2}=m_{1}\alpha,$ there exists $y\in\{k+1,k+2,\ldots\}$ such that $y\gamma_{2}=m_{1}\alpha+1.$ Since $(k-1)\gamma_{2}=m_{1}\alpha+1=y\gamma_{2}$ and $k-1\not=y,$ the transformation $\gamma_{2}$ is not injective. Moreover, $\left\vert \nb(\gamma_{2})\right\vert=0,\left\vert\{1,2,\ldots, n\}\gamma_{2}\right\vert=n,$ and $1\gamma_{2}=m_{1}\alpha+k-1\geq k> n.$ Hence, $\gamma_{2}\in\Lambda_{n}.$ 
By straightforward calculations, we obtain $A_{i}\gamma_{1}\gamma_{2}=\{m_{i}\alpha\}$ for all $1\leq i\leq l.$ Therefore, $\gamma_{1}\gamma_{2}=\alpha.$

Altogether, we have shown $F_{\mathbb{N}}\subseteq\Theta\Lambda_{n}.$ Since the converse inclusion is clear, we have $\Theta\Lambda_{n}=F_{\mathbb{N}}.$
\end{proof}
\end{proposition}
By the construction of $\gamma_{1}$ in Proposition \ref{Proposition4}, we observe that the only conditions for $\gamma_{1}$ are $M_{\alpha}=M_{\gamma_{1}}$ and $\min(\im \gamma_{1})\geq n.$ This gives us the following corollary.
\begin{corollary}\label{Corollary5}
Let $n\in\mathbb{N}$ and $\alpha\in F_{\mathbb{N}}.$ For $\gamma_{1}\in\Theta$ with $M_{\alpha}=M_{\gamma_{1}}$ and $\min(\im \gamma_{1})\geq n,$ there exists $\gamma_{2}\in\Lambda_{n}$ such that $\alpha=\gamma_{1}\gamma_{2}.$
\end{corollary}

As one can see, $F_{\mathbb{N}}$ is uncountable and thus, any generating set of $F_{\mathbb{N}}$ is uncountable. It appears the question whether a minimal generating set of $F_{\mathbb{N}}$ exists. The following constructions clarify that there are no minimal generating sets of $F_{\mathbb{N}},$ that is to say, we can get a smaller generating set (under the set inclusion) by excluding suitable elements from a given generating set. 

Let $\alpha\in F_{\mathbb{N}}^{\inf},$ $R_{\alpha}:=\{x\in\im\alpha:x\alpha^{-1}~\text{is not a convex set}\},$ and $Q_{\alpha}:=\{x\in\im\alpha:\left\vert x\alpha^{-1}\right\vert,\left\vert (x+1)\alpha^{-1}\right\vert\geq 3\}.$
Further, let $P:=\{\alpha\in F_{\mathbb{N}}^{\inf}:\left\vert \bigcup_{n>3}M_{\alpha}^{n}\right\vert,\left\vert R_{\alpha}\right\vert,\left\vert Q_{\alpha}\right\vert<\aleph_{0}\}.$ For $l\in\mathbb{N},$ let $$K_{l}:=\{\alpha\in P:\left\vert MS_{\alpha}^{l}\right\vert=\aleph_{0}~\text{and}~\left\vert MS_{\alpha}^{n}\right\vert<\aleph_{0}~\text{for all}~n<l\}.$$ Note that $\left\vert M_{\alpha}^{*}\right\vert=\aleph_{0}$ for all $\alpha\in K_{l}.$ Further, let $K_{\aleph_{0}}:=P\setminus\bigcup_{n\in\mathbb{N}}K_{n}.$ 
\begin{lemma}\label{1}
Let $\alpha\in F_{\mathbb{N}}^{\inf}$ with $\left\vert R_{\alpha}\right\vert<\aleph_{0}.$ Then there is $k\in\mathbb{N}$ such that $a\alpha\leq b\alpha$ for all $k\leq a<b.$
\begin{proof}
Since $\left\vert R_{\alpha}\right\vert<\aleph_{0},$ there is $k'\in\mathbb{N}$ such that $x\alpha^{-1}$ is a convex set for all $x\geq k'.$ Let $k=\min(k'\alpha^{-1})$ and let $a,b\in\mathbb{N}$ with $k\leq a<b.$ Assume that $a\alpha<k',$ i.e. $k<a.$ Then $\rank\alpha=\aleph_{0}$ implies that $\{a,a+1,\ldots\}\alpha$ is an infinite convex set containing $k',$ that is, there is $s>a$ with $s\alpha=k'.$ Thus, $k'\alpha^{-1}$ is not a convex set because $k<a<s,$ where $s,k\in k'\alpha^{-1}$ and $a\not\in k'\alpha^{-1},$ a contradiction. Hence, $k'\leq a\alpha.$ Assume $b\alpha<a\alpha.$ Then $\rank\alpha=\aleph_{0}$ implies that $\{b,b+1,\ldots\}\alpha$ is an infinite convex set containing $a\alpha,$ that is, there exists $t\in \mathbb{N}$ with $b<t$ and $t\alpha=a\alpha.$ This means that $(a\alpha)\alpha^{-1}$ is not a convex set since $a<b<t,$ where $a,t\in(a\alpha)\alpha^{-1}$ and $b\not\in(a\alpha)\alpha^{-1},$ a contradiction to $k'\leq a\alpha.$ Therefore, $a\alpha\leq b\alpha.$
\end{proof}
\end{lemma}
As a consequence of Lemma \ref{1}, we obtain that $\alpha|_{B}$ is injective for all $B\in MS_{\alpha}\cap C_{k}.$
\begin{lemma}\label{2}
Let $\alpha,\beta\in F_{\mathbb{N}}^{\inf}$ and let $x\in R_{\beta}$ with $x\beta^{-1}\subseteq\im\alpha.$ Then $x\in R_{\alpha\beta}.$ 
\begin{proof}
Assume $x\not\in R_{\alpha\beta}.$ This means that $x(\alpha\beta)^{-1}=x\beta^{-1}\alpha^{-1}$ is a convex set. Since $x\beta^{-1}\subseteq \im\alpha,$ we obtain that $x\beta^{-1}\alpha^{-1}\alpha=x\beta^{-1}$ is a convex set. That means $x\not\in R_{\beta},$ a contradiction. Hence, $x\in R_{\alpha\beta}.$
\end{proof}
\end{lemma}
\begin{lemma}\label{beforeideal}
Let $\beta\in F_{\mathbb{N}}^{\inf}$ and let $X\subseteq\mathbb{N}$ be such that $\left\vert X\right\vert=\aleph_{0}$ and $\left\vert X\beta\right\vert<\aleph_{0}.$ Then $\left\vert R_{\beta}\right\vert=\aleph_{0}.$ Moreover, $\left\vert R_{\alpha\beta}\right\vert=\aleph_{0}$ for all $\alpha\in F_{\mathbb{N}}^{\inf}.$
\begin{proof}
Assume $\left\vert R_{\beta}\right\vert<\aleph_{0}.$ By Lemma \ref{1}, there is $k\in\mathbb{N}$ with $a\alpha\leq b\alpha$ for all $k\leq a<b.$ Let $B=\{x\in X:x\geq k\}$ and $c=\max (B\beta).$ Then $\left\vert B\right\vert=\aleph_{0}.$ Let $t\in\mathbb{N}$ with $t\geq k.$ Since $\left\vert B\right\vert=\aleph_{0},$ there is $s\in B$ such that $t<s.$ Then $t\beta\leq s\beta\leq c.$ This implies that $\rank\beta\leq k+c<\aleph_{0},$ a contradiction. Hence, $\left\vert R_{\beta}\right\vert=\aleph_{0}$ and so $\left\vert\{x\in R_{\beta}:x\beta^{-1}\subseteq\im\alpha\}\right\vert=\aleph_{0}.$ Therefore, $\left\vert R_{\alpha\beta}\right\vert=\aleph_{0}$ by Lemma \ref{2}.
\end{proof}
\end{lemma}
\begin{proposition}\label{ideal}
$F_{\mathbb{N}}\setminus P$ is an ideal of $F_{\mathbb{N}}.$
\begin{proof}
Let $\alpha\in F_{\mathbb{N}}\setminus P$ and $\beta\in F_{\mathbb{N}}.$ If $\rank\alpha<\aleph_{0}$ or $\rank\beta<\aleph_{0},$ then we obtain that $\rank\alpha\beta,\rank\beta\alpha<\aleph_{0},$ that is, $\alpha\beta, \beta\alpha\in F_{\mathbb{N}}\setminus P.$ Suppose now $\rank\alpha=\rank\beta=\aleph_{0}.$ Since $\im\alpha$ and $\im\beta$ are convex sets, we have that $\rank\alpha\beta=\aleph_{0}$ and $\rank\beta\alpha=\aleph_{0},$ respectively. Let $M_{\beta}=\{B_{i}:i\in\mathbb{N}\}$ with $B_{i}<B_{i+1}$ for all $i\in\mathbb{N}.$\\

\noindent\textbf{Case 1:} $\left\vert R_{\alpha}\right\vert=\aleph_{0}.$ Suppose that $R_{\alpha}=\{x_{i}:i\in\mathbb{N}\}$ with $x_{i}<x_{i+1}$ for all $i\in\mathbb{N}.$ Let $r$ be the least $q\in\mathbb{N}$ with $\min(\im\beta)\leq\min (x_{q}\alpha^{-1})$ and let $E=\{x_{i}:i\geq r\}.$ Then $x\alpha^{-1}\subseteq\im\beta$ for all $x\in E.$ Therefore, Lemma \ref{2} implies that $x\in R_{\beta\alpha}$ and so $E\subseteq R_{\beta\alpha}.$ Hence, $\left\vert R_{\beta\alpha}\right\vert\geq \left\vert E\right\vert=\aleph_{0}.$ 

Suppose $\left\vert R_{\alpha\beta}\right\vert<\aleph_{0}.$ Then there is $k\in\mathbb{N}$ such that $x\beta^{-1}\alpha^{-1}$ is a convex set for all $x\geq k.$ Moreover, $\left\vert R_{\alpha}\beta\right\vert=\aleph_{0}.$ Otherwise $\left\vert R_{\alpha}\beta\right\vert<\aleph_{0}$ and so Lemma \ref{beforeideal} implies $\left\vert R_{\alpha\beta}\right\vert=\aleph_{0},$ a contradiction. Therefore, $\left\vert R_{\alpha}\beta\cap\{k,k+1,\ldots\}\right\vert=\aleph_{0}.$ Let $s$ be the least $q\in\mathbb{N}$ such that $\min(\im\alpha)<\min(x_{q}\beta\beta^{-1})$ and let $D=\{x_{i}:i\geq s\}\beta\cap\{k,k+1,\ldots\}.$ Let $x\in D.$ Then $x\beta^{-1}\alpha^{-1}$ is a convex set and $x\beta^{-1}\cap R_{\alpha}\not=\emptyset.$ Suppose that $x_{j}\in x\beta^{-1}\cap R_{\alpha}$ for some $j\in\mathbb{N}.$ If $x\beta^{-1}\cap R_{\alpha}=\{x_{j}\},$ then $x\beta^{-1}\alpha^{-1}=x_{j}\alpha^{-1}$ is not a convex set, a contradiction. Thus, $\left\vert x\beta^{-1}\cap\im\alpha\right\vert\geq 3.$ Since $x_{j}\alpha^{-1}$ is not a convex set, we obtain $\left\vert x_{j}\alpha^{-1}\right\vert\geq 2.$ Hence, $\left\vert x\beta^{-1}\alpha^{-1}\right\vert>3.$ Therefore, $\left\vert\bigcup_{n>3}M^{n}_{\alpha\beta}\right\vert\geq\left\vert D\right\vert=\aleph_{0}.$\\

\noindent\textbf{Case 2:} $\left\vert \bigcup_{n>3}M_{\alpha}^{n}\right\vert=\aleph_{0}$ and $\left\vert R_{\alpha}\right\vert<\aleph_{0}.$ Let $\bigcup_{n>3}M_{\alpha}^{n}=\{A_{i}:i\in\mathbb{N}\}$ with $A_{i}<A_{i+1}$ for all $i\in\mathbb{N}.$ Let $r$ be the least $q\in\mathbb{N}$ such that $\min(\im\beta)\leq\min(A_{q}).$ Then for $i\geq r,$ there is $m_{i}\in\mathbb{N}$ with $\big(\bigcup_{j=m_{i}}^{m_{i}+\left\vert A_{i}\right\vert-1}B_{j}\big)\beta\subseteq A_{i}.$ Hence, there is $D_{i}\in M_{\beta\alpha}$ with $\big(\bigcup_{j=m_{i}}^{m_{i}+\left\vert A_{i}\right\vert-1}B_{j}\big)\subseteq D_{i}.$ Then $\left\vert D_{i}\right\vert\geq \left\vert\bigcup_{j=m_{i}}^{m_{i}+\left\vert A_{i}\right\vert-1}B_{j}\right\vert\geq\left\vert A_{i}\right\vert>3.$ This shows that $\left\vert \bigcup_{n>3}M_{\beta\alpha}^{n}\right\vert\geq\left\vert\bigcup_{i\in\mathbb{N}}D_{i}\right\vert=\aleph_{0}.$ 

If $\left\vert\big(\bigcup_{i\in\mathbb{N}}A_{i}\big)\alpha\beta\right\vert=\aleph_{0},$ then we obtain $\left\vert \bigcup_{n>3}M_{\alpha\beta}^{n} \right\vert=\aleph_{0}.$ Suppose now that $\left\vert\big(\bigcup_{i\in\mathbb{N}}A_{i}\big)\alpha\beta\right\vert<\aleph_{0}.$ Assume $\left\vert\big(\bigcup_{i\in\mathbb{N}}A_{i}\big)\alpha\right\vert<\aleph_{0}.$ Let $X=\{\min (A_{i}):i\in\mathbb{N}\}.$ Then $\left\vert X\right\vert=\aleph_{0}$ and $\left\vert X\alpha\right\vert<\aleph_{0}.$ So, Lemma \ref{beforeideal} implies that $\left\vert R_{\alpha}\right\vert=\aleph_{0},$ a contradiction. Hence, $\left\vert\big(\bigcup_{i\in\mathbb{N}}A_{i}\big)\alpha\right\vert=\aleph_{0}.$ Then $\left\vert R_{\alpha\beta}\right\vert=\aleph_{0}$ by Lemma \ref{beforeideal}.\\

\noindent\textbf{Case 3:} $\left\vert Q_{\alpha}\right\vert=\aleph_{0}.$ Then $\left\vert Q_{\alpha}\cap\im\beta\alpha\right\vert=\aleph_{0}$ since $\rank\beta\alpha=\aleph_{0}.$ This implies that $\left\vert Q_{\beta\alpha}\right\vert=\aleph_{0}.$

Suppose that $\left\vert Q_{\alpha\beta}\right\vert,\left\vert R_{\alpha\beta}\right\vert<\aleph_{0}.$ Then $\left\vert Q_{\alpha}\beta\right\vert=\aleph_{0}.$ Otherwise $\left\vert Q_{\alpha}\beta\right\vert<\aleph_{0}$ and so Lemma \ref{beforeideal} implies $\left\vert R_{\alpha\beta}\right\vert=\aleph_{0},$ a contradiction. Let $Q_{\alpha}=\{x_{i}:i\in\mathbb{N}\}$ with $x_{i}<x_{i+1}$ for all $i\in\mathbb{N}.$ Since $\left\vert Q_{\alpha\beta}\right\vert,\left\vert R_{\alpha\beta}\right\vert<\aleph_{0},$ there is $k\in\mathbb{N}$ such that $x\beta^{-1}\alpha^{-1}$ is a convex set, and $\left\vert x\beta^{-1}\alpha^{-1}\right\vert<3$ or $\left\vert (x+1)\beta^{-1}\alpha^{-1}\right\vert<3$ for all $x\geq k.$ Then $\left\vert Q_{\alpha}\beta\cap\{k,k+1,\ldots\}\right\vert=\aleph_{0}$ since $\left\vert Q_{\alpha}\beta\right\vert=\aleph_{0}.$ Let $D=Q_{\alpha}\beta\cap\{k,k+1,\ldots\}$ and let $x\in D.$ Then there is $s\in Q_{\alpha}$ such that $s\beta=x.$ Since $s\in Q_{\alpha},$ we obtain that $\left\vert s\alpha^{-1}\right\vert,\left\vert (s+1)\alpha^{-1}\right\vert\geq 3.$ Assume that $(s+1)\beta\not=x.$ Then $(s+1)\beta=x+1.$ Otherwise, $(s+1)\beta=x-1$ and thus, there is $t>s+1$ with $t\beta=x.$ Hence, $x\beta^{-1}\alpha^{-1}$ is not a convex set with $x\geq k,$ a contradiction. Thus, $\left\vert x\beta^{-1}\alpha^{-1}\right\vert\geq\left\vert s\alpha^{-1}\right\vert\geq 3$ and $\left\vert (x+1)\beta^{-1}\alpha^{-1}\right\vert\geq\left\vert (s+1)\alpha^{-1}\right\vert\geq 3,$ a contradiction to $x\in D.$ Hence, $x=s\beta=(s+1)\beta,$ that is, $\left\vert x\beta^{-1}\alpha^{-1}\right\vert\geq\left\vert\{s,s+1\}\alpha^{-1}\right\vert\geq 6$ and so $x\in\bigcup_{n>3}M^{n}_{\alpha\beta}.$ Therefore, $\left\vert\bigcup_{n>3}M^{n}_{\alpha\beta}\right\vert\geq\left\vert D\right\vert=\aleph_{0}.$

For all three cases, we obtain that $\alpha\beta,\beta\alpha\not\in P.$ Therefore, we can conclude that $F_{\mathbb{N}}\setminus P$ is an ideal of $F_{\mathbb{N}}.$

\end{proof}
\end{proposition}
\begin{lemma}\label{Proposition1}
Let $\alpha\in K_{l}$ for some $l\in\mathbb{N}$ and let $G$ be a generating set of $F_{\mathbb{N}}.$ Then there are $\gamma_{1}\in K_{l_{1}}\cup K_{\aleph_{0}}$ and $\gamma_{2}\in K_{l_{2}}\cup K_{\aleph_{0}}$ for some $l_{1},l_{2}\in\mathbb{N}$ with $l_{1},l_{2}>l$ such that $\alpha=\gamma_{1}\gamma_{2}$ and $\gamma_{1},\gamma_{2}\in\langle G\setminus\{\alpha\}\rangle.$
\begin{proof}
Since $\alpha\in K_{l},$ we have $\left\vert M_{\alpha}^{*}\right\vert=\aleph_{0}.$ Suppose that $M_{\alpha}^{*}=\{B_{i}: i\in\mathbb{N}\}$ with $B_{i}<B_{i+1}$ for all $i\in\mathbb{N}.$ Then we define $\gamma_{1}:\mathbb{N}\to\mathbb{N}$ by $\im\gamma_{1}=\mathbb{N},\gamma_{1}\in\Theta,$ and $M_{\gamma_{1}}^{*}=\{B_{i}:i\in2\mathbb{N}\}.$ Note that $\gamma_{1}$ is well defined by these three conditions. Moreover, we define $\gamma_{2}:\mathbb{N}\to\mathbb{N}$ by $x\gamma_{2}:=(\min(x\gamma_{1}^{-1}))\alpha$ for all $x\in\mathbb{N}.$ By the definitions of $\gamma_{1}$ and $\gamma_{2},$ it is clear that $\gamma_{1}\gamma_{2}=\alpha$ and that there exist $l_{1},l_{2}>l$ such that $\gamma_{1}\in K_{l_{1}}\cup K_{\aleph_{0}}$ and $\gamma_{2}\in K_{l_{2}}\cup K_{\aleph_{0}}.$ Hence, for $i\in\{1,2\},$ there is $k_{i}\in\mathbb{N}$ satisfying the following properties:
\begin{itemize}
\item[(a1)] $\left\vert A\right\vert\geq l_{i}>l$ for all $A\in MS_{\gamma_{i}}\cap C_{k_{i}};$
\item[(a2)] $\left\vert A \right\vert=3$ for all $A\in M_{\gamma_{i}}^{*}\cap C_{k_{i}};$
\item[(a3)] $\left\vert x\gamma_{i}^{-1}\right\vert<3$ or $\left\vert (x+1)\gamma_{i}^{-1}\right\vert<3$ for all $x\geq k_{i}\gamma_{i};$ 
\item[(a4)] $x\gamma_{i}^{-1}$ is a convex set for all $x\geq k_{i}\gamma_{i}$
\end{itemize}
because $\left\vert\bigcup_{n=1}^{l_{i}-1} MS^{n}_{\gamma_{i}}\right\vert<\aleph_{0}$ with $l_{i}>l,$ $\left\vert\bigcup_{n>3} M^{n}_{\gamma_{i}}\right\vert<\aleph_{0},\left\vert Q_{\gamma_{i}}\right\vert<\aleph_{0},$ and $\left\vert R_{\gamma_{i}}\right\vert<\aleph_{0},$ respectively. It is a consequence of (a4) that $a\gamma_{i}\leq b\gamma_{i}$ for all $k_{i}\leq a<b,$ which we will use without further mention.
Since $\alpha\in K_{l},$ there is $k\in\mathbb{N}$ satisfying the following properties:
\begin{itemize}
\item[(b1)] $\left\vert MS_{\alpha}^{l}\cap C_{k}\right\vert=\aleph_{0};$
\item[(b2)] $\left\vert A\right\vert=3$ for all $A\in M_{\alpha}^{*}\cap C_{k}$
\end{itemize}
because $\left\vert MS^{l}_{\alpha}\right\vert=\aleph_{0}$ and $\left\vert\bigcup_{n>3} M^{n}_{\alpha}\right\vert<\aleph_{0},$ respectively.
Since $\langle G\rangle=F_{\mathbb{N}}$ and $\gamma_{1},\gamma_{2}\in P,$ there are $\mu_{1},\mu_{2},\ldots,\mu_{m_{1}},\eta_{1},\eta_{2},\ldots,\eta_{m_{2}}\in G\cap P$ such that $\gamma_{1}=\mu_{1}\mu_{2}\cdots\mu_{m_{1}}$ and $\gamma_{2}=\eta_{1}\eta_{2}\cdots\eta_{m_{2}}$ for some $m_{1},m_{2}\in\mathbb{N}.$ By (a1) and (b1), it is clear that $\mu_{1}\not=\alpha$ and $\eta_{1}\not=\alpha.$

Assume that $\mu_{j}=\alpha$ for some $j\in\{2,3,\ldots,m_{1}\}.$ Let $MS_{\alpha}^{l,k}=\{A\in MS_{\alpha}^{l}: \{k\}<A\}=\{A_{i}:i\in\mathbb{N}\}$ with $A_{i}<A_{i+1}$ for all $i\in\mathbb{N}.$ Let $x\in\mathbb{N}$ be such that $x>k_{1}+3$ and $x\delta_{1}\in\{\min(A):A\in MS_{\alpha}^{l,k}\setminus\{A_{1}\}\}.$ Let $\delta_{1}=\mu_{1}\mu_{2}\cdots\mu_{j-1}.$ Further, let $\delta_{2}=\mu_{j+1}\mu_{j+2}\cdots\mu_{m_{1}}$ if $j<m_{1}$ and let $\delta_{2}=\identity_{\mathbb{N}}$ if $j=m_{1}.$ Note that $\identity_{\mathbb{N}}\in P.$ Then $x\delta_{1}=\min(A_{r})$ for some $r\geq 2$ and so $A_{r}=\{x\delta_{1},x\delta_{1}+1,\ldots,x\delta_{1}+l-1\}.$ So, (b2) implies that $B_{1}=\{x\delta_{1}-3,x\delta_{1}-2,x\delta_{1}-1\},B_{2}=\{x\delta_{1}+l,x\delta_{1}+l+1,x\delta_{1}+l+2\}\in M_{\alpha}.$ Note that $k<x-3.$ 
 
Since $\{x-3,x-2,x-1,x\}\delta_{1}$ is a convex set containing $x\delta_{1},$ we get that $\{x-3,x-2,x-1\}\delta_{1}\subseteq B_{1}$ and so $\{x-3,x-2,x-1\}\subseteq (x-1)\delta_{1}\alpha\delta_{2}(\delta_{1}\alpha\delta_{2})^{-1}.$ We obtain the equality $\{x-3,x-2,x-1\}= (x-1)\delta_{1}\alpha\delta_{2}(\delta_{1}\alpha\delta_{2})^{-1}$ by (a2). Let $D=\{x,x+1,\ldots,x+l_{1}-1\}.$ Note that $z\gamma_{1}\gamma_{1}^{-1}$ is a convex set for all $z\in D.$ By (a3), we can conclude that $\left\vert x\delta_{1}\alpha\delta_{2}(\delta_{1}\alpha\delta_{2})^{-1}\right\vert=\left\vert x\gamma_{1}\gamma_{1}^{-1}\right\vert=1.$ Let $A=\{X\in M_{\gamma_{1}}^{*}: X\subseteq D\setminus\{x\}\}.$ Assume that $A\not=\emptyset.$ Then there is $E\in A$ with $E\leq X$ for all $X\in A.$ Then $\{x,x+1,\ldots,\min(E)-1\}\in\bigcup_{n=1}^{l_{1}-1} MS^{n}_{\delta_{1}\alpha\delta_{2}},$ a contradiction. This implies that $\delta_{1}|_{D}$ is injective with $z\delta_{1}=x\delta_{1}+z-x$ for all $z\in D.$ Since $l_{1}>l,$ we have $x+l\in D$ with $(x+l)\delta_{1}\alpha\alpha^{-1}=(x\delta_{1}+l)\alpha\alpha^{-1}=B_{2}.$ Then $(x+l)\gamma_{1}\gamma_{1}^{-1}=(x+l)\delta_{1}\alpha\delta_{2}(\delta_{1}\alpha\delta_{2})^{-1}=(x\delta_{1}+l)\alpha\delta_{2}\delta_{2}^{-1}\alpha^{-1}\delta_{1}^{-1}\supseteq(x\delta_{1}+l)\alpha\alpha^{-1}\delta_{1}^{-1}=B_{2}\delta_{1}^{-1}.$ Therefore, $\left\vert (x+l)\gamma_{1}\gamma_{1}^{-1}\right\vert\geq\left\vert B_{2}\delta_{1}^{-1}\right\vert\geq\left\vert B_{2}\right\vert=3,$ a contradiction. Therefore, we conclude that $\mu_{j}\not=\alpha$ for all $j\in\{1,2,\ldots,m_{1}\}.$ Similarly, we can show that $\eta_{j}\not=\alpha$ for all $j\in\{1,2,\ldots,m_{2}\}.$ So, $\gamma_{1},\gamma_{2}\in\langle G\setminus\{\alpha\}\rangle.$
 \end{proof}
\end{lemma}
In particular, Lemma \ref{Proposition1} shows that $G$ has no common elements to $K_{l}$ for all $l\in\mathbb{N},$ whenever $G$ is a minimal generating set of $F_{\mathbb{N}}.$ The main result of this section states that there are no minimal generating sets of $F_{\mathbb{N}}.$ If such a one existed, it would have the following necessary condition.
\begin{lemma}\label{Corollary2}
If $G$ is a minimal generating set of $F_{\mathbb{N}}$, then $G\cap K_{n}=\emptyset$ for all $n\in\mathbb{N}.$ Moreover, $G\cap P\subseteq K_{\aleph_{0}}.$
\begin{proof}
Assume $G\cap K_{l}\not=\emptyset$ for some $l\in\mathbb{N}.$ Then there exists $\alpha\in G\cap K_{l.}$ By Lemma \ref{Proposition1}, there are $\gamma_{1},\gamma_{2}\in\langle G\setminus\{\alpha\}\rangle$ with $\alpha=\gamma_{1}\gamma_{2},$ that is, $\alpha\in\langle G\setminus \{\alpha\}\rangle.$ Since $\langle G\rangle=F_{\mathbb{N}},$ we obtain $\langle G\setminus\{\alpha\}\rangle=F_{\mathbb{N}}.$ It contradicts to the assumption that $G$ is a minimal generating set of $F_{\mathbb{N}}.$ Therefore, $G\cap K_{n}=\emptyset$ for all $n\in\mathbb{N}.$ Together with $P=\big(\bigcup_{n\in\mathbb{N}}K_{n}\big)\cup K_{\aleph_{0}},$ we obtain that $G\cap P=G\cap\big(\big(\bigcup_{n\in\mathbb{N}}K_{n}\big)\cup K_{\aleph_{0}}\big)=G\cap K_{\aleph_{0}}\subseteq K_{\aleph_{0}}.$ 
\end{proof}
\end{lemma}
\begin{theorem}\label{Proposition3}
There are no minimal generating sets of $F_{\mathbb{N}}.$
\begin{proof}
Assume that there is a minimal generating set $G$ of $F_{\mathbb{N}}.$ By Lemma \ref{Corollary2}, we have $G\cap K_{n}=\emptyset$ for all $n\in\mathbb{N}.$ Now, we define $\alpha:\mathbb{N}\to\mathbb{N}$ by
$$x\alpha:=\begin{cases}
2n-1~&\text{if}~x=4n-3~\text{for}~n\in\mathbb{N};\\
2n&\text{if}~x\in\{4n-2,4n-1,4n\}~\text{for}~n\in\mathbb{N}.
\end{cases}$$
Then $M_{\alpha}^{*}=\{\{4n-2,4n-1,4n\}: n\in\mathbb{N}\}.$ It is clear that $\alpha\in P$ since $R_{\alpha}=Q_{\alpha}=\bigcup_{n>3}M^{n}_{\alpha}=\emptyset.$ Since $\alpha\in P$ and $\langle G\rangle=F_{\mathbb{N}},$ Lemma \ref{Corollary2} implies that $\alpha=\gamma_{1}\gamma_{2}\cdots\gamma_{l}$ for some $\gamma_{1},\gamma_{2},\ldots,\gamma_{l}\in G\cap P\subseteq K_{\aleph_{0}}$ and for some $l\in\mathbb{N}.$ Let $\gamma_{0}=\identity_{\mathbb{N}}$ be the identity mapping on $\mathbb{N}$ and let $i\in\{1,2,\ldots,l\}.$ Since $\alpha=\gamma_{1}\gamma_{2}\cdots\gamma_{l},$ we obtain the following properties:
\begin{itemize}
\item[(a1)] $a\gamma_{i}\leq b\gamma_{i}$ for all $1\gamma_{0}\gamma_{1}\cdots\gamma_{i-1}\leq a<b;$
\item[(a2)] $\left\vert B\right\vert=3$ for all $B\in M_{\gamma_{i}}^{*}\cap C_{1\gamma_{0}\gamma_{1}\cdots\gamma_{i-1}}$
\end{itemize}
because $R_{\alpha}=\emptyset$ and $M_{\alpha}^{*}=M_{\alpha}^{3},$ respectively. Moreover, (a1) provides 
\begin{itemize}
\item[(a3)] $\gamma_{i}|_{A}$ is injective for all $A\in MS_{\gamma_{i}}\cap C_{1\gamma_{0}\gamma_{1}\cdots\gamma_{i-1}}.$
\end{itemize} 
Let $a_{l}=2$ and $a_{l-j}=2a_{l-j+1}+3$ for all $j\in\mathbb{N}\setminus\{l,l+1,\ldots\}.$ Since $\gamma_{i}\in K_{\aleph_{0}},$ there exists $m_{i}\in\mathbb{N}$ such that $\left\vert C\right\vert\geq a_{i}$ for all $C\in MS_{\gamma_{i}}\cap C_{m_{i}}.$ Let $m^{*}=\max\{1\gamma_{1},1\gamma_{1}\gamma_{2},\linebreak\ldots,1\gamma_{1}\gamma_{2}\cdots\gamma_{l-1},m_{1},m_{2},\ldots,m_{l}\}$ and let $y\in\mathbb{N}$ be such that $\{m^{*}\}<\{y,y\gamma_{1},y\gamma_{1}\gamma_{2},\ldots,\linebreak y\gamma_{1}\gamma_{2}\cdots\gamma_{l-1}\}.$ Further, let $D_{1}\in MS_{\gamma_{1}}\cap C_{y}$ and let $x=\min(D_{1}).$ Since $m^{*}<y\leq x,$ we obtain that $\left\vert D_{1}\right\vert\geq a_{1}$ and $\gamma_{1}|_{D_{1}}$ is injective by (a3). Let $j\in\{2,3,\ldots,l\}.$ Then $m^{*}<y\leq x$ and (a1) imply that $m^{*}\leq y\gamma_{1}\gamma_{2}\cdots\gamma_{j-1}\leq x\gamma_{1}\gamma_{2}\cdots\gamma_{j-1}.$ Since $a_{j-1}=2a_{j}+3$ and $m^{*}\leq x\gamma_{1}\gamma_{2}\cdots\gamma_{j-1},$ the properties (a2) and (a3) provide that there is a convex set $D_{j}\subseteq D_{j-1}\gamma_{j-1}\cap E_{j}$ for some $E_{j}\in MS_{\gamma_{j}}$ such that $\left\vert D_{j}\right\vert=a_{j}$ and $\gamma_{j}|_{D_{j}}$ is injective. Let $D=D_{l}\gamma_{l-1}^{-1}\gamma_{l-2}^{-1}\cdots\gamma_{1}^{-1}.$ Since $D\gamma_{0}\gamma_{1}\cdots\gamma_{r-1}\subseteq D_{r},\gamma_{r}|_{D_{r}}$ is injective, and $D_{r}\gamma_{r}\gamma_{r}^{-1}=D_{r}$ for all $1\leq r\leq l,$ we obtain that $\left\vert D\right\vert=\left\vert D_{l}\right\vert=a_{l}=2.$ Then there is $D'\in MS_{\gamma_{1}\gamma_{2}\cdots\gamma_{l}}$ with $D\subseteq D'.$ Thus, $\left\vert D'\right\vert \geq \left\vert D\right\vert=2,$ a contradiction to $\alpha=\gamma_{1}\gamma_{2}\cdots\gamma_{l}$ with $MS_{\alpha}=MS_{\alpha}^{1}.$

\end{proof}
\end{theorem}
Although a minimal generating set of the uncountable semigroup $F_{\mathbb{N}}$ does not exist, there is an uncountable subsemigroup of $F_{\mathbb{N}}$ having such one. Let $A\subseteq\mathbb{N}$ and let $\alpha_{A}\in\Theta$ be such that  $\im\alpha_{A}=\mathbb{N}$ and $\left\vert x\alpha_{A}^{-1} \right\vert=3$ if $x\in A$ and $\left\vert x\alpha_{A}^{-1} \right\vert=5$ otherwise. Obviously, $\alpha_{A}$ is well defined. Further, let $Q:=\{\alpha_{A}: A\subseteq\mathbb{N}\}.$ Then $\left\vert Q\right\vert=2^{\aleph_{0}},$ which means that $Q$ is uncountable. For $A,B\subseteq\mathbb{N},$ it is easy to verify that $\left\vert M_{\alpha_{A}\alpha_{B}}^{m}\right\vert>0$ for some $m\geq 9,$ that is, $\alpha_{A}\alpha_{B}\not\in Q.$ This shows that $Q$ is a minimal generating set of the semigroup generated by $Q.$ In other words, the uncountable subsemigroup $\langle Q\rangle$ of $F_{\mathbb{N}}$ has a minimal generating set. 
 
\section{Infinite Decreasing Chains of Generating Sets of $F_{\mathbb{N}}$}
The previous section shows that there are no minimal generating sets of $F_{\mathbb{N}}.$ Obviously, $F_{\mathbb{N}}$ itself is the maximum generating set. Both facts provide that $F_{\mathbb{N}}$ must have infinite decreasing chains of generating sets of $F_{\mathbb{N}}.$ In this section, we will provide such two chains.

Let $\Inj(F_{\mathbb{N}})$ be the set of all injective transformations in $F_{\mathbb{N}}$ and let $\xi$ be the transformation on $\mathbb{N}$ defined by $x\xi:=x+2$ for all $x\in\mathbb{N}.$ Thus, $\xi^{n}\in\Inj(F_{\mathbb{N}})$ with $1\xi^{n}=2n+1$ for all $n\in\mathbb{N}.$ Let $\mathcal{B}:=\{\alpha\in F_{\mathbb{N}}:\left\vert\nb(\alpha)\right\vert =2, \collapse(\alpha)=3,~\text{and}~\im\alpha=\mathbb{N}\}.$ For $n\in\mathbb{N},$ there is exactly one $\beta\in\mathcal{B}$ with $\min(\nb(\beta))= n.$ This transformation will be denoted by $\beta_{n}.$ Let $n\in\mathbb{N}.$ We put $\mathcal{B}_{n}:=\{\beta_{i}:i\geq n\}$ and $n\chi:=\begin{cases}
1&\text{if }n~\text{is odd};\\
2&\text{otherwise}.
\end{cases}$ \\Further, we define transformations $\lambda_{n}$ and $\delta_{n}$ as follows:
$$x\lambda_{n}:=\begin{cases}
n-x+1~~&\text{if }x\in\{1,2,\ldots,n\};\\
x-n+1&\text{otherwise}
\end{cases}$$
and
$$x\delta_{n}:=\begin{cases}
n\chi~~&\text{if }x\in\{1,2,\ldots,n\};\\
n\chi+x-n&\text{otherwise.}
\end{cases}$$
It is easy to check that $\delta_{n}\in F_{\mathbb{N}}.$ But $\lambda_{n}\in F_{\mathbb{N}},$ whenever $n$ is odd. In this case, we observe that $\left\vert \nb(\lambda_{n})\right\vert=0,\left\vert\{1,2,\ldots, n\}\lambda_{n}\right\vert=n,$ and $1\lambda_{n}=n.$ If $n\not=1,$ then $(n-1)\lambda_{n}=2=(n+1)\lambda_{n},$ that is, $\collapse(\lambda_{n})>0$ and so $\lambda_{n}\in\Lambda_{n}.$
\begin{lemma}\label{Proposition6}
Let $n\in\mathbb{N}.$ Then $\delta_{m}\in\langle \mathcal{B}_{n}\cup \Lambda_{n}\cup \{\xi\} \rangle$ for all $m\in\mathbb{N}.$
\begin{proof}
Let $m\in\mathbb{N},m_{1}=\max\{m,n\},$ and $m_{2}=2m_{1}+1.$ Then we can calculate that
$$\delta_{m}=\begin{cases}
\xi\beta_{1}~~&\text{if }m=n=1;\\
\xi^{m_{1}}\beta_{m_{2}-2}\lambda_{m_{2}-2}&\text{if }m=1,n>1;\\
\xi^{m_{1}}\beta_{m_{2}}^{k_{1}}\lambda_{m_{2}}&\text{if }m=2k_{1}+1~\text{for some}~k_{1}\in\mathbb{N};\\
\xi^{m_{1}}\beta_{m_{2}-1}^{k_{2}}\lambda_{m_{2}-2}&\text{if }m=2k_{2}~\text{for some}~k_{2}\in\mathbb{N}.
\end{cases}$$
Clearly, $\beta_{1}\in\mathcal{B}_{1}.$ If $n+m>2,$ then $m_{2}-2>n,$ which implies that $\beta_{m_{2}-2},\beta_{m_{2}-1},\beta_{m_{2}}\in\mathcal{B}_{n}$ and $\lambda_{m_{2}-2}, \lambda_{m_{2}}\in\Lambda_{n}.$ Altogether, we obtain $\delta_{m}\in\langle \mathcal{B}_{n}\cup\Lambda_{n}\cup\{\xi\}\rangle.$
\end{proof}
\end{lemma}
\noindent Let $n\in\mathbb{N}.$ We define a transformation $\alpha_{n}$ on $\mathbb{N}$ by $x\alpha_{n}:=x$ if $x\in\mathbb{N}\setminus\{n,n+1,\ldots\}$ and $x\alpha_{n}:=n$ otherwise. It is clear that $\alpha_{n}\in F_{\mathbb{N}}.$ Then we put $\mathcal{A}_{n}:=\{\alpha_{i}:i\geq n\}.$ Further, let $$\Delta:=\{\alpha\in F_{\mathbb{N}}:\left\vert M_{\alpha}^{*}\right\vert=\aleph_{0}\}$$ and $\Delta_{n}:=\Delta\cap\Omega_{n}=\{\alpha\in F_{\mathbb{N}}:1\alpha\geq n,\left\vert \{1,2,\ldots,n\}\alpha\right\vert=n,~\text{and}~\left\vert M_{\alpha}^{*}\right\vert=\aleph_{0}\}.$ 
\begin{lemma}\label{Proposition7}
Let $\alpha\in F_{\mathbb{N}}\setminus\Delta.$ Then $\alpha\in\langle\mathcal{A}_{n}\cup\mathcal{B}_{n}\cup\Lambda_{n}\cup\{\xi\} \rangle$ for all $n\in\mathbb{N}.$ 
\begin{proof}
Since $\alpha\in F_{\mathbb{N}}\setminus\Delta,$ we have $\left\vert M_{\alpha}^{*}\right\vert<\aleph_{0}.$ Let $n\in\mathbb{N}$ and let $k_{1}\in\mathbb{N}\setminus\{1,2,\ldots,n\}$ be odd. Further, let $k'=\frac{1}{2}(k_{1}-1).$\\
\noindent Suppose that $\left\vert M_{\alpha}^{*}\right\vert=0.$ Then $\left\vert \nb(\alpha)\right\vert=0.$ Thus, $x$ and $x\alpha$ have the same parity for all $x\in\mathbb{N}.$ We define $\gamma:\mathbb{N}\to\mathbb{N}$ by
$$x\gamma:=\begin{cases}
1\alpha+k_{1}-x~~&\text{if }x\in\{1,2,\ldots, k_{1}-1\};\\
(x-k_{1}+1)\alpha&\text{otherwise}.
\end{cases}$$
Then $\left\vert \nb(\gamma)\right\vert=0,\collapse(\gamma)>0,1\gamma=1\alpha+k_{1}-1>n,$ and $\left\vert \{1,2,\ldots,n\}\gamma\right\vert=n,$ that is, $\gamma\in\Lambda_{n}.$ So, we obtain $\alpha=\xi^{k_{1}^{'}}\gamma\in\langle \Lambda_{n}\cup\{\xi\}\rangle.$\\

\noindent Suppose now $M_{\alpha}^{*}=\{A_{i}:1\leq i\leq m\}$ for some $m\in\mathbb{N}$ with $A_{i}< A_{j}$ for all $1\leq i<j\leq m.$ It follows $\left\vert A_{i}\right\vert<\aleph_{0}$ for all $i\in\mathbb{N}\setminus\{m,m+1,\ldots\}.$ Let 
$$p_{i}=\min(A_{i})~\text{for all}~i\in\{1,2,\ldots,m\}$$
and
$$m_{i}=\max(A_{i})~\text{for all}~i\in\mathbb{N}\setminus\{m,m+1,\ldots\}.$$ 
Further, let $k_{i+1}=k_{i}+p_{i+1}-m_{i}$ for all $i\in\mathbb{N}\setminus\{m,m+1,\ldots\}.$ If $i\in\mathbb{N}\setminus\{1,m,m+1,\ldots\},$ then $1\not\in A_{i}$ and $\left\vert A_{i}\right\vert=2l_{i}+1$ for some $l_{i}\in\mathbb{N}.$ If $1\not\in A_{1},$ then we suppose $\left\vert A_{1}\right\vert=2l_{1}+1$ for some $l_{1}\in\mathbb{N}.$ In the case $\left\vert A_{m}\right\vert<\aleph_{0},$ we obtain that $\left\vert A_{m}\right\vert=2l_{m}+1$ for some $l_{m}\in\mathbb{N}.$ We define transformations $\gamma_{1},\gamma_{2},\ldots,\gamma_{m}$ on $\mathbb{N}$ as follows:
$$\gamma_{1}:=\begin{cases}
\delta_{m_{1}}\xi^{k_{1}^{'}}~~~&\text{if }1\in A_{1};\\
\xi^{k_{1}^{'}}\beta_{k_{1}+p_{1}-1}^{l_{1}}&\text{otherwise,}
\end{cases}$$
for $i\in\mathbb{N}\setminus\{1,m,m+1,\ldots\},$ we put
$$\gamma_{i}:=\begin{cases}
\beta_{k_{i}}^{l_{i}}~~~&\text{if }1\in A_{1}~\text{and}~m_{1}~\text{is odd;}\\
\beta_{k_{i}+1}^{l_{i}}&\text{if }1\in A_{1}~\text{and}~m_{1}~\text{is even;}~~~~~~~~~~~~~~\\
\beta_{k_{i}+p_{1}-1}^{l_{i}}~&\text{if }1\not\in A_{1},
\end{cases}$$
and
$$\gamma_{m}:=\begin{cases}
\beta_{k_{m}}^{l_{m}}~~&\text{if }1\in A_{1}, m_{1}~\text{is odd, and}~\left\vert A_{m}\right\vert<\aleph_{0};\\
\alpha_{k_{m}}&\text{if }1\in A_{1}, m_{1}~\text{is odd, and}~\left\vert A_{m}\right\vert=\aleph_{0};\\
\beta_{k_{m}+1}^{l_{m}}&\text{if }1\in A_{1}, m_{1}~\text{is even, and}~\left\vert A_{m}\right\vert<\aleph_{0};\\
\alpha_{k_{m}+1}&\text{if }1\in A_{1}, m_{1}~\text{is even, and}~\left\vert A_{m}\right\vert=\aleph_{0};\\
\beta_{k_{m}+p_{1}-1}^{l_{m}}&\text{if }1\not\in A_{1}~\text{and}~\left\vert A_{m}\right\vert<\aleph_{0};\\
\alpha_{k_{m}+p_{1}-1}&\text{if }1\not\in A_{1}~\text{and}~\left\vert A_{m}\right\vert=\aleph_{0}.
\end{cases}$$
Let $\alpha^{*}=\gamma_{1}\gamma_{2}\cdots\gamma_{m}.$ By straightforward calculations, we obtain that $\alpha^{*}\in\Theta,$ $M_{\alpha}=M_{\alpha^{*}},$ and $1\alpha^{*}\geq k_{1}> n.$ Then Corollary \ref{Corollary5} implies that there exists $\alpha^{'}\in\Lambda_{n}$ with $\alpha=\alpha^{*}\alpha^{'}.$ By the definition of $\gamma_{1}$ and Lemma \ref{Proposition6}, we get $\gamma_{1}\in\langle \mathcal{B}_{n}\cup\Lambda_{n}\cup\{\xi\}\rangle.$ For $i\in\{2,3,\ldots,m\},$ we obtain that $\gamma_{i}\in\langle \mathcal{A}_{n}\cup\mathcal{B}_{n}\rangle$ since $k_{i}> n.$ Therefore, $\alpha=\alpha^{*}\alpha^{'}\in\langle \mathcal{A}_{n}\cup\mathcal{B}_{n}\cup\Lambda_{n}\cup\{\xi\}\rangle.$
\end{proof}
\end{lemma}
Both previous lemmas lead to the definition of an infinite decreasing chain $\{H_{n}:n\in\mathbb{N}\}$ of generating sets of $F_{\mathbb{N}},$ where $H_{n}:=\mathcal{A}_{n}\cup\mathcal{B}_{n}\cup\Lambda_{n}\cup\Delta_{n}\cup\{\xi\}.$ It is worth mentioning  that the intersection of the $H_{i}$'s gives the singleton set $\{\xi\},$ which is not a generating set of $F_{\mathbb{N}}.$ It is easy to verify that $\xi\not\in\langle\mathcal{A}_{n}\cup\mathcal{B}_{n}\cup\Lambda_{n}\cup\Delta_{n}\rangle.$ Therefore, the relative rank of $F_{\mathbb{N}}$ modulo $\mathcal{A}_{n}\cup\mathcal{B}_{n}\cup\Lambda_{n}\cup\Delta_{n}$ is one.
\begin{theorem}\label{Proposition8}
$\langle H_{n}\rangle=F_{\mathbb{N}}$ for all $n\in\mathbb{N}.$
\begin{proof}
Let $n\in\mathbb{N}.$ It is a consequence of Lemma \ref{Proposition7} that $$\langle\mathcal{A}_{n}\cup\mathcal{B}_{n}\cup\Lambda_{n}\cup\Delta\cup\{\xi\}\rangle=F_{\mathbb{N}}.$$ In order to show $\langle H_{n}\rangle=F_{\mathbb{N}},$ it is enough to prove $\Delta\setminus\Delta_{n}\subseteq \langle H_{n}\rangle.$ Let $\alpha\in\Delta\setminus\Delta_{n}.$ Then $\left\vert M_{\alpha}^{*}\right\vert=\aleph_{0}$ and so $\left\vert M_{\alpha}\right\vert=\aleph_{0}.$ Suppose that $M_{\alpha}=\{A_{i}:i\in\mathbb{N}\}$ with $A_{i}<A_{i+1}$ for all $i\in\mathbb{N}.$ Let $p_{i}=\min(A_{i})$ for all $i\in\mathbb{N}$ and let $k_{1}\in\mathbb{N}$ be odd such that $k_{1}>n.$\\ 

\noindent\textbf{Case 1:} $\left\vert\{1,2,\ldots,n\}\alpha\right\vert=n.$ We define $\gamma:\mathbb{N}\to\mathbb{N}$ by $x\gamma:=k_{1}+i-1~\text{for all }x\in A_{i}, i\in\mathbb{N}.$ It is obvious that $\gamma\in\Theta, M_{\gamma}^{*}=M_{\alpha}^{*},1\gamma=k_{1}>n,$ and $\left\vert\{1,2,\ldots,n\}\gamma\right\vert=n.$ This means $\gamma\in\Delta_{n}.$ Moreover, Corollary \ref{Corollary5} implies that there exists $\gamma'\in\Lambda_{n}$ with $\gamma\gamma'=\alpha.$ Therefore, $\alpha\in\langle H_{n}\rangle.$\\

\noindent\textbf{Case 2:} $\left\vert \{1,2,\ldots,n\}\alpha\right\vert<n.$ Let $s$ be the smallest natural number $r$ such that $n<p_{r}$ and $A_{r}\in M_{\alpha}^{*}.$ Then we define $\gamma_{0}:\mathbb{N}\to\mathbb{N}$ by
$$x\gamma_{0}:=\begin{cases}
k_{1}+x-1~~~~~~~&\text{if }x\in\{1,2,\ldots,p_{s}-1\};\\
k_{1}+p_{s}+i-2&\text{if }x\in A_{s+i-1}~\text{for }i\in\mathbb{N}.
\end{cases}$$
Note that $\gamma_{0}\in\Delta_{n}$ since $1\gamma_{0}=k_{1}>n,\left\vert\{1,2,\ldots,n\}\gamma_{0}\right\vert=n,$ and $\left\vert M_{\gamma_{0}}^{*}\right\vert=\left\vert M_{\alpha}^{*}\right\vert-s=\aleph_{0}.$ If $s=\min\{i\in\mathbb{N}: A_{i}\in M_{\alpha}^{*}\},$ then $M_{\gamma_{0}}=M_{\alpha}$ and so we put $\beta:=\gamma_{0}.$ Suppose $s>\min\{i\in\mathbb{N}: A_{i}\in M_{\alpha}^{*}\}.$ Let $\{C\in M_{\alpha}^{*}:C<A_{s}\}=\{B_{i}:1\leq i\leq m\}$ for some $m\in\mathbb{N}$ with $B_{i}<B_{j}$ for all $1\leq i<j\leq m.$ For $i\in\mathbb{N}\setminus\{1,m+1,m+2,\ldots\},$ there is $l_{i}\in\mathbb{N}$ with $\left\vert B_{i}\right\vert=2l_{i}+1.$ Moreover, there is $l_{1}\in\mathbb{N}$ with $\left\vert B_{1}\right\vert=2l_{1}+1$ or $\left\vert B_{1}\right\vert=2l_{1},$ depending on the parity of $\left\vert B_{1}\right\vert.$ Let $q_{i}=\min(B_{i})$ and $m_{i}=\max(B_{i})$ for all $i\in\{1,2,\ldots,m\}.$ Further, let $k_{j+1}=k_{j}+q_{j+1}-m_{j}$ for all $j\in\mathbb{N}\setminus\{m,m+1,\ldots\}.$ For $i\in\{1,2,\ldots,m\},$ we define $\gamma_{i}:\mathbb{N}\to\mathbb{N}$ as follows:
$$\gamma_{i}:=\begin{cases}
\beta_{k_{i}}^{l_{i}}~~&\text{if }1\in B_{1}~\text{and}~\left\vert B_{1}\right\vert~\text{is odd;}\\
\beta_{k_{i}-1}^{l_{i}}&\text{if }1\in B_{1}~\text{and}~\left\vert B_{1}\right\vert~\text{is even;}\\
\beta_{k_{i}+q_{i}-1}^{l_{i}}&\text{if }1\not\in B_{1}.
\end{cases}$$ 
In this case, we put $\beta:=\gamma_{0}\gamma_{1}\gamma_{2}\cdots\gamma_{m}.$ By straightforward calculations, we obtain that $\beta\in\Theta,M_{\beta}=M_{\alpha},$ and $1\beta\geq k_{1}-1\geq n.$  Then Corollary \ref{Corollary5} implies that there exists $\beta^{'}\in\Lambda_{n}$ such that $\beta\beta^{'}=\alpha.$ Therefore, $\alpha=\beta\beta^{'}\in\langle H_{n}\rangle.$
\end{proof}
\end{theorem}
It is easy to see that $\Omega_{n+1}\subsetneq\Omega_{n},\mathcal{A}_{n+1}\subsetneq\mathcal{A}_{n},$ and $\mathcal{B}_{n+1}\subsetneq\mathcal{B}_{n}$ for all $n\in\mathbb{N}.$ Therefore, we can conclude that $\{H_{n}:n\in\mathbb{N}\}$ is an infinite decreasing chain of generating sets of $F_{\mathbb{N}}.$

Recall that $F_{\mathbb{N}}=\Theta\Lambda_{n}$ for any $n\in\mathbb{N},$ where $\Theta$ is a subsemigroup of $F_{\mathbb{N}}.$ This means that we can generate any element in $F_{\mathbb{N}}$ by elements from $\Theta$ and $\Lambda_{n}.$ Now, let $$\Gamma:=\{\alpha\in \Theta: \rank\alpha=\aleph_{0}~\text{and there exists}~b\in\im\alpha~\text{with}~\left\vert b\alpha^{-1}\right\vert\geq 3\}.$$
We will generate the elements in $F_{\mathbb{N}}$ by elements from the proper subsemigroup $\Gamma$ of $F_{\mathbb{N}},\Lambda_{n},$ and the additional transformation $\xi,$ for any $n\in\mathbb{N}.$ Moreover, $\Lambda_{n}$ is covered by the semigroup $\Lambda.$

\begin{proposition}\label{Proposition9}
$\Lambda$ and $\Gamma$ are subsemigroups of $F_{\mathbb{N}}.$
\begin{proof}
Let $\alpha,\beta\in\Lambda.$ Then $\left\vert\nb(\alpha)\right\vert=\left\vert\nb(\beta)\right\vert=0$ and $\collapse(\alpha),\collapse(\beta)>0.$ This means $M_{\alpha}^{*}=M_{\beta}^{*}=\emptyset.$ Assume $\left\vert M_{\alpha\beta}^{*}\right\vert>0.$ Then there exists $D\in M_{\alpha\beta}^{*},$ that is, $\left\vert D\right\vert>1$ and $\left\vert D\alpha\beta\right\vert=1.$ Since $D$ is a convex set and $\left\vert D\right\vert>1,$ there is $a\in\mathbb{N}$ such that $\{a,a+1\}\subseteq D.$ Since $\left\vert\nb(\alpha)\right\vert=0,$ we obtain that $a\alpha=b$ and $(a+1)\alpha=c$ for some $b,c\in\mathbb{N}$ such that $\left\vert b-c\right\vert=1.$ Since $\left\vert\{b,c\}\beta\right\vert=\left\vert \{a,a+1\}\alpha\beta\right\vert\leq\left\vert D\alpha\beta\right\vert=1$ and $\left\vert b-c\right\vert=1,$ we obtain $\left\vert\nb(\beta)\right\vert\not=0,$ a contradiction. Therefore, $M_{\alpha\beta}^{*}=\emptyset,$ that is, $\left\vert \nb(\alpha\beta)\right\vert=0.$ Together with $0<\collapse(\alpha)\leq\collapse(\alpha\beta),$ we obtain that $\alpha\beta\in\Lambda.$\\

Now, let $\alpha,\beta\in\Gamma.$ Then $\alpha,\beta\in\Theta$ and $\rank\alpha=\rank\beta=\aleph_{0}.$ It is clear that $\rank\alpha\beta=\aleph_{0}$ and $\alpha\beta\in\Theta.$ Furthermore, there is $a\in\mathbb{N}$ with $\left\vert a\alpha^{-1}\right\vert\geq 3.$ Then $\left\vert a\beta(\alpha\beta)^{-1}\right\vert=\left\vert a\beta\beta^{-1}\alpha^{-1}\right\vert\geq\left\vert a\alpha^{-1}\right\vert\geq 3.$ Altogether, we conclude that $\alpha\beta\in\Gamma.$
\end{proof} 
\end{proposition}
We are going to establish a second infinite decreasing chain of generating sets of $F_{\mathbb{N}}$, which are subsets of the union of the three semigroups $\{\xi\},\Lambda,$ and $\Gamma.$ Let $n\in\mathbb{N}$ and let $G_{n}$ be the set of all $\alpha\in F_{\mathbb{N}}$ satisfying at least one of the following three properties:
\begin{itemize}
\item[(g1)] $\alpha=\xi;$
\item[(g2)] $\alpha\in \Lambda_{n};$
\item[(g3)] $\alpha\in\Theta_{n}$ such that $\left\vert M_{\alpha}^{*}\right\vert\in\{1,\aleph_{0}\}$ and $M_{\alpha}^{*}=M_{\alpha}^{3}.$ 
\end{itemize}
Clearly, $G_{n}\subseteq \Gamma\cup\Lambda_{n}\cup\{\xi\}.$
\begin{theorem}\label{Theorem10}
$\langle G_{n}\rangle=F_{\mathbb{N}}$ for all $n\in\mathbb{N}.$
\begin{proof}
Let $n\in\mathbb{N}.$ By the definition of $G_{n}$, we have $\Lambda_{n}\cup\{\xi\}\subseteq G_{n}.$ We will show that $\mathcal{A}_{n},\mathcal{B}_{n},\Delta_{n}\subseteq \langle G_{n}\rangle.$\\

Let $\alpha\in\mathcal{A}_{n}.$ Then $\alpha=\alpha_{k}$ for some $k\geq n,$ and $x\alpha=x$ if $x\in\mathbb{N}\setminus\{k,k+1,\ldots\}$ and $x\alpha=k$ otherwise. Let $l$ be the least even natural number $r$ such that $r>k.$ We define transformations $\gamma_{1}$ and $\gamma_{2}$ on $\mathbb{N}$ as follows:
$$x\gamma_{1}:=\begin{cases}
l+x~~&\text{if }x\in\mathbb{N}\setminus\{k,k+1,\ldots\};\\
l+k&\text{if }x\in\{k,k+2,k+4,\ldots\};\\
l+k+1&\text{if }x\in\{k+1,k+3,k+5,\ldots\}
\end{cases}$$
and
$$~~~~~x\gamma_{2}:=\begin{cases}
l+x~~&\text{if }x\in\{1,2,\ldots,l+k-1\};\\
2l+k&\text{if }x\in\{l+k,l+k+1,l+k+2\};\\
l+x-2&\text{if }x\in\mathbb{N}\setminus\{1,2,\ldots,l+k+2\}.
\end{cases}$$
Then $\gamma_{1}\in\Lambda_{n}$ and $\gamma_{2}$ satisfies (g3). By straightforward calculations, we obtain $\gamma_{1}\gamma_{2}\lambda_{2l+1}\\=\alpha.$ Since $1\lambda_{2l+1}=2l+1>n,$ we have $\lambda_{2l+1}\in\Lambda_{n}.$ This shows $\mathcal{A}_{n}\subseteq\langle G_{n}\rangle.$\\

Let $\alpha\in\mathcal{B}_{n}.$ Then $\alpha=\beta_{k}$ for some $k\geq n,$ that is,
$$x\alpha=\begin{cases}
x~~&\text{if }x\in\mathbb{N}\setminus\{k,k+1,\ldots\};\\
k&\text{if }x\in\{k,k+1,k+2\};\\
x-2&\text{if }x\in\mathbb{N}\setminus\{1,2,\ldots,k+2\}. 
\end{cases}$$
Let $l$ be again the least even natural number $r$ such that $r>k$ and define $\gamma:\mathbb{N}\to\mathbb{N}$ by $x\gamma:=x\alpha+l$ for all $x\in\mathbb{N}.$ Then $\gamma$ satisfies (g3). It is easy to see that $\gamma\lambda_{l+1}=\alpha.$ Since $1\lambda_{l+1}=l+1>n,$ we obtain $\lambda_{l+1}\in\Lambda_{n},$ that is, $\mathcal{B}_{n}\subseteq\langle G_{n}\rangle.$\\

Let $\alpha\in\Delta_{n}.$ Then $1\alpha\geq n, \left\vert\{1,2,\ldots,n\}\alpha\right\vert=n,$ and $\left\vert M_{\alpha}^{*}\right\vert=\aleph_{0}.$ Suppose $M_{\alpha}^{*}=\{A_{i}:i\in\mathbb{N}\}$ with $A_{i}<A_{i+1}$ for all $i\in\mathbb{N}.$ It follows that $\left\vert A_{i}\right\vert<\aleph_{0}$ for all $i\in\mathbb{N}.$ For $i\in\mathbb{N},$ let $p_{i}=\min(A_{i})$ and $l_{i}=\left\vert A_{i}\right\vert.$ Let $l$ be now the least even natural number $r$ such that $r>1\alpha.$ Further, let $k_{2}=l+p_{2}$ and $k_{i}=l+p_{i}-\Sigma_{j=2}^{i-1}(l_{j}-3)$ for all $i\in\mathbb{N}\setminus\{1,2\}.$ Note that if $l_{1}$ is even, then $p_{1}=1.$ Put $c=1$ if $l_{1}$ is even and $c=0$ otherwise. We define transformations $\gamma_{1},\gamma_{2},$ and $\gamma_{3}$ on $\mathbb{N}$ as follows:
$$x\gamma_{1}:=\begin{cases}
x~~~~~~&\text{if }x\in\{1,2,\ldots,p_{2}-1\};\\
k_{i}&\text{if }x\in\{p_{i},p_{i}+2,\ldots,p_{i}+l_{i}-3\};\\
k_{i}+1&\text{if }x\in\{p_{i}+1,p_{i}+3,\ldots,p_{i}+l_{i}-2\};\\
k_{i}+2&\text{if }x=p_{i}+l_{i}-1;\\
l+x-\Sigma_{j=1}^{i}(l_{j}-3)&\text{if }x\in\{p_{i}+l_{i},p_{i}+l_{i}+1,\ldots,p_{i+1}-1\},~~~~~~~~~~~~~~~~~~~
\end{cases}$$
$$x\gamma_{2}:=\begin{cases}
l+x+l_{1}-3+c~~&\text{if }x\in\{1,2,\ldots,l+p_{1}-1-c\};\\
2l+p_{1}+l_{1}-3&\text{if }x\in\{l+p_{1}-c,l+p_{1}+2-c,\ldots,l+p_{1}+l_{1}-3\};\\
2l+p_{1}+l_{1}-2&\text{if }x\in\{l+p_{1}+1-c,l+p_{1}+3-c,\ldots,l+p_{1}+l_{1}-2\};\\
l+x&\text{if }x\in\{l+p_{1}+l_{1}-1,l+p_{1}+l_{1},\ldots\},
\end{cases}$$
and
$$x\gamma_{3}:=\begin{cases}
l+x&\text{if }x\in\{1,2,\ldots,2l+p_{1}+l_{1}-4\};\\
3l+p_{1}+l_{1}-3&\text{if }x\in\{2l+p_{1}+l_{1}-3,2l+p_{1}+l_{1}-2,2l+p_{1}+l_{1}-1\};\\
l+x-2&\text{if }x\in\{2l+p_{1}+l_{1},2l+p_{1}+l_{1}+1,\ldots,l+k_{2}-1\};\\
2l+k_{i}-2(i-1)&\text{if }x\in\{l+k_{i},l+k_{i}+1,l+k_{i}+2\};\\
l+x-2i&\text{if }x\in\{l+k_{i}+3,l+k_{i}+4,\ldots,l+k_{i+1}-1\}
\end{cases}$$
for all $i\in\mathbb{N}\setminus\{1\}.$ It is easy to verify that $\gamma_{1},\gamma_{2}\in\Lambda_{n}$ and $\gamma_{3}$ satisfies (g3). By straightforward calculations, we obtain that $\gamma_{1}\gamma_{2}\gamma_{3}\in\Theta,M_{\gamma_{1}\gamma_{2}\gamma_{3}}=M_{\alpha},$ and $1\gamma_{1}\gamma_{2}\gamma_{3}\geq 2l+l_{1}-2\geq l>n.$ Then Corollary \ref{Corollary5} implies that there exists $\gamma_{4}\in\Lambda_{n}$ such that $\gamma_{1}\gamma_{2}\gamma_{3}\gamma_{4}=\alpha.$ Therefore, $\Delta_{n}\subseteq \langle G_{n}\rangle.$\\

Altogether, we have shown $H_{n}=\mathcal{A}_{n}\cup\mathcal{B}_{n}\cup\Lambda_{n}\cup\Delta_{n}\cup\{\xi\}\subseteq \langle G_{n}\rangle.$ By Proposition \ref{Proposition8}, we obtain $\langle G_{n}\rangle=F_{\mathbb{N}}.$
\end{proof}
\end{theorem}

Let $n\in\mathbb{N}.$ Since $\Omega_{n+1}\subsetneq\Omega_{n},$ we can conclude that $G_{n+1}\subsetneq G_{n}.$ This shows that $\{G_{n}:n\in\mathbb{N}\}$ is an infinite decreasing chain of generating sets of $F_{\mathbb{N}}.$ Moreover, $\bigcap_{n\in\mathbb{N}} G_{n}=\{\xi\}$ because any transformation $\alpha\in F_{\mathbb{N}}\setminus\{\xi\}$ is not in $G_{1\alpha+1}.$ In other words, the relative rank of $F_{\mathbb{N}}$ modulo $G_{n}$ is one.

\vskip1cm \centerline{ \textbf{Acknowledgments}}

The first author would like to express her thanks to the Development and Promotion of Science and Technology Talents Project and the Department of Mathematics, Faculty of Science, Khon Kaen University. 

\end{document}